\documentclass[12pt]{amsart}
\usepackage{amscd,amssymb}
\newtheorem{thm}{Theorem}[section]

\newtheorem{prop}[thm]{Proposition}
\theoremstyle{definition}

\theoremstyle{remark}
\newtheorem{remark}[thm]{Remark}

\newtheorem{example}[thm]{Example}

\newtheorem{Question}[thm]{Question}
\numberwithin{equation}{section}

\def\R {{\Bbb R}}

\def\Homeo{{\mathrm{Homeo}}\,}

\def\QED{\nobreak\quad\ifmmode\roman{Q.E.D.}\else{\rm Q.E.D.}\fi}

\def\Exp{\operatorname{Exp}}
\def\sbs{\subset}
\def\sS{{\mathcal{S}}}

\begin{document}

\title{On universal minimal compact $G$-spaces}
\author[V. Uspenskij]{Vladimir Uspenskij}
\address{Department of Mathematics,
321 Morton Hall, Ohio University, Athens, Ohio 45701, USA}
\email{uspensk@bing.math.ohiou.edu}

\date{June 10, 2000}

\begin{abstract} 
For every topological group $G$ one can define the universal minimal compact
$G$-space $X=M_G$ characterized by the following properties: 
(1) $X$ has no proper
closed $G$-invariant subsets; (2) for every compact $G$-space $Y$ there exists
a $G$-map $X\to Y$. If $G$ is the group of all orientation-preserving
homeomorphisms of the circle $S^1$, then $M_G$ can be identified with $S^1$
(V. Pestov). We show that the circle cannot be replaced by the Hilbert cube
or a compact manifold of dimension $>1$. This answers a question of V. Pestov. Moreover, we prove that for every topological group $G$ 
the action of $G$ on $M_G$ is not 3-transitive.
\end{abstract}
\maketitle

\setcounter{tocdepth}{1}

\section{Introduction}
With every topological group $G$ one can associate the {\em universal
minimal compact $G$-space} $M_G$. To define this object, recall some basic
definitions. A {\em $G$-space} is a topological space $X$ with a continuous
action of $G$, that is, a map $G\times X\to X$ satisfying $g(hx)=(gh)x$
and $1x=x$ ($g,h\in G$, $x\in X$). A $G$-space $X$ is {\em minimal} if it has
no proper $G$-invariant closed subsets or, equivalently, if the orbit
$Gx$ is dense in $X$ for every $x\in X$. A map $f:X\to Y$ between two
$G$-spaces is {\em $G$-equivariant}, or a {\em $G$-map} for short, if
$f(gx)=gf(x)$ for every $g\in G$ and $x\in X$. 

All {\em maps} are assumed to be continuous, and `compact' includes
`Hausdorff'. The universal minimal compact $G$-space $M_G$ is characterized
by the following property: $M_G$ is a minimal compact $G$-space, and 
for every compact minimal $G$-space $X$ there exists a $G$-map of $M_G$
onto $X$. Since Zorn's lemma implies that every compact $G$-space has 
a minimal compact $G$-subspace, it follows that for every compact $G$-space
$X$, minimal or not, there exist a $G$-map of $M_G$ to $X$.

The existence of $M_G$ is easy: consider the product of a representative
family of compact minimal $G$-spaces, and take any minimal closed $G$-subspace
of this product for $M_G$. It is also true that $M_G$ is unique, in the sense
that any two universal minimal 
compact $G$-spaces are isomorphic \cite{Aus}. 
For the reader's
convenience, we give a proof of this fact in the Appendix.

If $G$ is locally compact, the action of $G$ on $M_G$ is free
\cite{Veech} (see also \cite{P3}, Theorem 3.1.1), 
that is, if $g\ne1$, then $gx\ne x$ for every $x\in M_G$. 
On the other hand, $M_G$ is a singleton 
for many naturally arising non-locally compact groups $G$. 
This property of $G$ is equivalent
to the following {\em fixed point on compacta (f.p.c.) property}: 
every compact $G$-space has a $G$-fixed point. (A point
$x$ is $G$-fixed if $gx=x$ for all $g\in G$.) For example, if $H$ is
a Hilbert space, the group $U(H)$ of all unitary operators on $H$, equipped
with the pointwise convergence topology, has the f.p.c. property (Gromov --
Milman); another example of a group with this property, due to Pestov, 
is $H_+(\R)$, the
group of all orientation-preserving self-homeomorphisms of the real line.
We refer the reader to beautiful papers by V.~Pestov
\cite{P1,P2,P3} on this subject.

Let $S^1$ be a circle, and let $G=H_+(S^1)$ be the group of all 
orientation-preserving self-homeomorphisms of $S^1$. Then $M_G$ can be
identified with $S^1$ \cite{P1}, Theorem 6.6. 
The question arises whether 
a similar assertion holds for the Hilbert cube $Q$. This question is 
due to V.~Pestov, who writes in  \cite{P1}, Concluding Remarks,
that his theorem ``tends to suggest that the Hilbert cube $I^\omega$
might serve as the universal minimal flow for the group $\Homeo(I^\omega)$". 
In other words, let $G=H(Q)$ be the group
of all self-homeomorphisms of $Q=I^\omega$, equipped with the compact-open
topology. Are $M_G$ and $Q$ isomorphic as $G$-spaces? 

The aim of the present paper is to answer this question in the negative. 
Let us say that the action of a group $G$ on a $G$-space $X$ is {\em
3-transitive} if $|X|\ge3$ and
for any triples $(a_1,a_2,a_3)$ and $(b_1,b_2,b_3)$ of
distinct points in $X$ there exists $g\in G$ such that $ga_i=b_i$, $i=1,2,3$.

\begin{thm}
For every topological group $G$ the action of $G$ on the universal minimal
compact $G$-space $M_G$ is not 3-transitive.
\label{main}
\end{thm}

Since the action of $H(Q)$ on $Q$ is 3-transitive, it follows that 
$M_G\ne Q$ for $G=H(Q)$. Similarly, if $K$ is compact and $G$ is
a 3-transitive group of homeomorphisms of $K$, then $M_G\ne K$. 
This remark
applies, for example, if $K$ is a manifold of dimension $>1$ or a Menger
manifold and $G=H(K)$.

\begin{Question}
Let $G=H(Q)$. Is $M_G$ metrizable?
\label{q1}
\end{Question}

A similar question
can be asked when $Q$ is replaced by a compact manifold or a Menger manifold.

Let $P$ be the pseudoarc (= the unique hereditarily indecomposable chainable
continuum) and $G=H(P)$. The action of $G$ on $P$ is transitive but not
2-transitive, and the following question remains open:

\begin{Question}
Let $P$ be the pseudoarc and $G=H(P)$. Can $M_G$ be identified with~$P$?
\label{pseudo}
\end{Question}

\section{Proof of the main theorem}
The proof of Theorem~\ref{main} depends on the consideration of the space
of maximal chains of closed sets. For a compact space $K$ let $\Exp K$ be
the (compact) space of all non-empty closed subsets of $K$, equipped with
the Vietoris topology. A subset $C\sbs \Exp K$ is a {\em chain} if for any
$E,F\in C$ either $E\sbs F$ or $F\sbs E$. 
If $C\sbs \Exp K$ is
a chain, so is the closure of $C$. It follows that every maximal chain is
a closed subset of $\Exp K$ and hence an element of $\Exp\Exp K$. Let 
$\Phi\sbs\Exp\Exp K$ be the space of all maximal chains. Then $\Phi$ is
closed in $\Exp\Exp K$ and hence compact. Let us sketch a proof. 
Clearly the closure of $\Phi$ consists of chains. 
Assume $C\in\Exp\Exp K$
is a non-maximal chain. 
We construct a neighbourhood $\frak W$ of $C$ in $\Exp\Exp K$ which 
is disjoint from $\Phi$. One the following cases holds:
(1) the first member of $C$ has more
than one point, or (2) the last member of $C$ is not $K$, or 
(3) the chain $C$ contains ``big gaps": there
are $F_1,F_2\in C$ such that $|F_2\setminus F_1|\ge2$ and 
for every 
$F\in C$ either $F\sbs F_1$ or $F_2\sbs F$.
For example, consider the third case (the first two cases are simpler).
Find open sets $U$, $V_1$, $V_2$ in $K$ with pairwise disjoint closures
such that $F_1\sbs U$ and $F_2$ meets both $V_1$ and $V_2$. Let $\frak W=
\{D\in \Exp\Exp K: \text{every member of $D$ either is contained in $U$ or
meets both $V_1$ and $V_2$\,}\}$. Then $\frak W$ is a neighbourhood of $C$
which does not meet $\Phi$. Indeed, suppose $D\in\frak W\cap \Phi$. 
Let $E_1$ be the largest member of $D$ which is contained in $\bar U$.
Let $E_2$ be the smallest
member of $D$ which meets both $\bar V_1$ and $\bar V_2$. 
For every $E\in D$
we have either $E\sbs E_1$ or $E_2\sbs E$, and $|E_2\setminus E_1|\ge 2$.
Pick a point $p\in E_2\setminus E_1$. The set $E_1\cup\{p\}$ is comparable
with every member of $D$ but is not a member of $D$. This contradicts the
maximality of $D$. We have proved that $\Phi$ is compact.
 
Suppose $G$ is a topological group and $K$ is a compact $G$-space.
Then the natural action of $G$ on $\Exp K$ is continuous, hence $\Exp K$
is a compact $G$-space, and so is $\Exp\Exp K$. Since the closed set 
$\Phi\sbs \Exp\Exp K$ is $G$-invariant, $\Phi$ is a compact $G$-space, too.

\begin{prop}
Let $G$ be a topological group. Pick $p\in M_G$, and let $H=\{g\in G:gp=p\}$
be the stabilizer of $p$. There exists a maximal chain $C$ of closed subsets
of $M_G$ such that $C$ is $H$-invariant: if $F\in C$ and $g\in H$, then 
$gF\in C$.
\label{prop}
\end{prop}

Note that members of an $H$-invariant chain need not be $H$-invariant.

\begin{proof}
Every compact $G$-space $X$ has an $H$-invariant point. Indeed, there
exists a $G$-map $f:M_G\to X$, and since $p$ is $H$-invariant, so is 
$f(p)\in X$.

Let $\Phi\sbs \Exp\Exp M_G$ 
be the compact space of all maximal chains of closed subsets of $M_G$.
We saw that $\Phi$ is a compact $G$-space. Thus $\Phi$ has an $H$-invariant
point.
\end{proof}

Theorem~\ref{main} follows from Proposition~\ref{prop}:

{\it Proof of Theorem~\ref{main}.} Assume that the action of $G$ on $X=M_G$
is 3-transitive. Pick $p\in X$, and let $H=\{g\in G:gp=p\}$. According to
Proposition~\ref{prop}, there exists an $H$-invariant maximal chain $C$
of closed subsets of $X$. 
The smallest member of $C$ is an $H$-invariant
singleton. Since $G$ is 2-transitive on $X$, the only $H$-invariant
singleton is $\{p\}$. Thus $\{p\}\in C$, and all members of $C$ contain $p$. 
Our definition of 3-transitivity implies that $|X|\ge3$. Thus there exists
$F\in C$ such that $F\ne \{p\}$ and $F\ne X$. Pick $a\in F\setminus\{p\}$ 
and $b\in X\setminus F$. The points $p,a,b$ are all distinct.
Since $G$ is 3-transitive on $X$,
there exists $g\in G$ such that $gp=p$, $ga=b$ and $gb=a$. 
Since $a\in F$ and $b\notin F$, we have $b=ga\in gF$ and $a=gb\notin gF$.
Thus $a\in F\setminus gF$ and $b\in gF\setminus F$, so $F$ and $gF$ are not
comparable.
On the other hand, the equality $gp=p$ means that $g\in H$. 
Since $C$ is $H$-invariant and $F\in C$, we have $gF\in C$. 
Hence $F$ and $gF$ must be
comparable, being members of the chain $C$. We have arrived at a contradiction.
\qed

\smallskip
\begin{example}
Consider the group $G=H_+(S^1)$ of all 
orientation-preserving self-homeomorphisms of the circle $S^1$. 
According to Pestov's result cited above, $M_G=S^1$. This example shows
that the action of $G$ on $M_G$ can be 2-transitive. Pick $p\in S^1$, and let
$H\sbs G$ be the stabilizer of $p$. Proposition~\ref{prop} implies that 
there must exist $H$-invariant maximal chains of closed subsets of $S^1$.
It is easy to see that there are precisely two such chains. They consist
of the singleton $\{p\}$, the whole circle and of all arcs that either
``start at $p$" or ``end at $p$", respectively.
\end{example}

\begin{remark}
Let $P$ be the pseudoarc, and let $G=H(P)$. Pick a point $p\in P$, and
let $H\sbs G$ be the stabilizer of $p$. Then there exists an $H$-invariant
maximal chain $C$ of closed subsets of $P$. Namely, let $C$ be the collection
of all subcontinua $F\sbs P$ such that $p\in F$. Since any two subcontinua
of $P$ are either disjoint or comparable, it follows that $C$ is a chain.
The chain $C$ can be shown to be maximal, and it is obvious that $C$ is
$H$-invariant. Thus Proposition~\ref{prop} does not contradict the conjecture
that $M_G=P$. This observation motivates our question~\ref{pseudo}.
\end{remark}

\section{Appendix: Uniqueness of $M_G$}
We sketch a proof of the uniqueness of $M_G$ up to a $G$-isomorphism.

Let $G$ be a topological group. The {\em greatest ambit} 
$X=\sS(G)$ for $G$
is a compact $G$-space with a distinguished point $e$ such that for 
every pointed compact $G$-space $(Y, e')$ there exists a unique 
$G$-morphism $f:X\to Y$ such that $f(e)=e'$. The greatest ambit is
defined uniquely up to a $G$-isomorphism preserving distinguished points. 
We can take for $\sS(G)$ the Samuel
compactification of $G$ equipped with the right uniformity, which is
the compactification of $G$ corresponding to the algebra of all bounded
right uniformly continuous functions. The distinguished point is the unity
of $G$. See \cite{P1, P2, P3} for more details. 

The greatest ambit $X$ has a natural structure of a left-topological semigroup. 
This means that there is an associative multiplication 
$(x,y)\mapsto xy$ on $X$ (extending the original multiplication on $G$)
such that for every $y\in X$ the self-map $x\mapsto xy$ of $X$ is continuous. 
Let $x, y\in X$.
There is a unique $G$-map $r_y:X\to X$ such that $r_y(e)=y$. Define
$xy=r_y(x)$. If $y$ is fixed, the map $x\mapsto xy$ is equal to $r_y$ and
hence is continuous. 
If $y,z\in X$, the self-maps $r_zr_y$ and $r_{yz}$ of $X$
are equal, since both are $G$-maps sending $e$ to $yz=r_z(y)$. 
This means that the 
multiplication on $X$ is associative. The distinguished element $e\in X$
is the unity of $X$: we have $ex=r_x(e)=x$ and $xe=r_e(x)=x$. 
If $g\in G$ and $x\in X$, the expression $gx$ can be
understood in two ways: in the sense of the exterior action of $G$ on $X$
and as a product in $X$; these two meanings agree.
If $f:X\to X$ is a $G$-self-map and $a=f(e)$, then 
$f(x)=f(xe)=xf(e)=xa=r_a(x)$ for all $x\in X$. 
Thus the semigroup of all $G$-self-maps of $X$ coincides
with the semigroup $\{r_y:y\in X\}$ of all right multiplications.

A subset $I\sbs X$ is a {\em left ideal} if $XI\sbs I$. Closed 
$G$-subspaces of $X$ are the same as closed left ideals of $X$.
An element $x$ of a semigroup is an {\it idempotent\/} if $x^2=x$. 
Every closed $G$-subspace of $X$, being a left ideal, is moreover 
a left-topological compact semigroup and hence contains an idempotent,
according to the following fundamental result of R. Ellis
(see \cite{Rup}, Proposition~2.1 or \cite{BJM}, Theorem~3.11):

\begin{thm}
Every non-empty compact left-topological semigroup $K$ contains an idempotent.
\label{idemp}
\end{thm}

\begin{proof}
Zorn's lemma implies that there exists a minimal element $Y$ in the set of all
closed non-empty subsemigroups of $K$. Fix $a\in Y$. We claim that $a^2=a$
(and hence $Y$ is a singleton). The set $Ya$, being a closed subsemigroup
of $Y$, is equal to $Y$. It follows that the closed subsemigroup $Z=\{x\in Y:
xa=a\}$ is non-empty. Hence $Z=Y$ and $xa=a$ for every $x\in Y$. In particular,
$a^2=a$.
\end{proof}

Let $M$ be a minimal closed left ideal of $X$. We have just proved that 
there is an idempotent $p\in M$. Since $Xp$ is a closed left ideal 
contained in $M$, we have $Xp=M$. Thus the $G$-map $r_p:X\to M$ defined by
$r_p(x)=xp$ is a retraction of $X$ onto $M$. In particular, $xp=x$ for 
every $x\in M$.

\begin{prop}
Every $G$-map $f:M\to M$ has the form $f(x)=xy$ for some $y\in M$
\label{pr2}
\end{prop}

\begin{proof}
The composition $h=fr_p:X\to M$ is a $G$-map of $X$ into itself, hence
it has the form $h=r_y$, where $y=h(e)\in M$. Since $r_p\restriction M
=\text{Id}$, we have $f=h\restriction M=r_y\restriction M$.
\end{proof}

\begin{prop}
Every $G$-map $f:M\to M$ is bijective.
\label{pr3}
\end{prop}

\begin{proof}
According to Proposition~\ref{pr2}, there is $a\in M$ such that $f(x)=xa$
for all $x\in M$. Since $Ma$ is a compact $G$-space contained in $M$, we
have $Ma=M$ by the minimality of $M$. Thus there exists $b\in M$ such that
$ba=p$. Let $g:M\to M$ be the $G$-map defined by $g(x)=xb$. Then 
$fg(x)=xba=xp=x$ for every $x\in M$, therefore
$fg=1$ (the identity map of $M$). We have proved that in the semigroup $S$
of all $G$-self-maps of $M$, every element has a right inverse. Hence $S$
is a group. (Alternatively, we first deduce from the equality $fg=1$ that
all elements of $S$ are surjective and then, applying this to $g$,
we see that $f$ is also injective.)
\end{proof}

We are now in a position to prove that {\em every universal compact minimal
$G$-space is isomorphic to $M$}. First note that the minimal compact $G$-space
$M$ is itself universal: if $Y$ is any compact $G$-space, there exists
a $G$-map of the greatest ambit $X$ to $Y$, and its restriction to $M$
is a $G$-map of $M$ to $Y$. Now let $M'$ be another universal compact minimal
$G$-space. There exist $G$-maps $f:M\to M'$ and $g:M'\to M$. Since $M'$ is
minimal, $f$ is surjective. On the other hand, in virtue of 
Proposition~\ref{pr3} the composition $gf:M\to M$ is bijective. It follows
that $f$ is injective and hence a $G$-isomorphism between $M$ and $M'$.

\end{document}